\documentclass[12pt]{article}
\usepackage{latexsym,amsmath,amssymb}
\usepackage[usenames]{color}

\newfont{\bb}{msbm10 at 12pt}
\newfont{\bbt}{msbm10 at 10pt}
\def\tr{\hbox{\bbt R}}

\def\r{\hbox{\bb R}}

\def\c{\hbox{\bb C}}

\newcommand{\beq}{\begin{equation}}

\newcommand{\eeq}{\end{equation}}

\newcommand{\norm}[1]{\left\Vert #1 \right\Vert}

\newcommand{\set}[1]{\left\{#1\right\}}

\newcommand{\eps}{\varepsilon}

\newcommand{\To}{\longrightarrow }

\usepackage[latin1]{inputenc}
\usepackage {amsmath}
\usepackage {amsthm}
\usepackage{times}
\usepackage{amscd}
\usepackage{epsf}

\numberwithin{equation} {section}

\begin{document}

\theoremstyle{plain}\newtheorem{lemma}{Lemma}[section]
\theoremstyle{plain}\newtheorem{definition}{Definition}[section]
\theoremstyle{plain}\newtheorem{theorem}{Theorem}[section]
\theoremstyle{plain}\newtheorem{proposition}{Proposition}[section]
\theoremstyle{plain}\newtheorem{remark}{Remark}[section]
\theoremstyle{plain}\newtheorem{corollary}{Corollary}[section]

\begin{center}
\rule{15cm}{1.5pt} \vspace{.6cm}

{\Large \bf A Colding-Minicozzi stability inequality \\[2mm] and its applications}
\vspace{0.4cm}

\vspace{0.5cm}

{\large Jos$\acute{\text{e}}$ M. Espinar$\,^\dag$\footnote{The author is partially
supported by Spanish MEC-FEDER Grant MTM2007-65249, and Regional J. Andalucia Grants
P06-FQM-01642 and FQM325}, Harold Rosenberg$\,^\ddag$}\\
\vspace{0.3cm} \rule{15cm}{1.5pt}
\end{center}

\vspace{.5cm}

\noindent $\mbox{}^\dag$ Departamento de Geometría y Topología, Universidad de
Granada, 18071 Granada, Spain; e-mail: jespinar@ugr.es\vspace{0.2cm}

\noindent $\mbox{}^\ddag$ Instituto de Matematica Pura y Aplicada, 110 Estrada Dona
Castorina, Rio de Janeiro 22460-320, Brazil; e-mail: rosen@impa.br

\vspace{.3cm}

\begin{abstract}
We consider operators $L$ acting on functions on a Riemannian surface, $\Sigma$, of
the form
$$L = \Delta + V -a K .$$

Here $\Delta$ is the Laplacian of $\Sigma$, $V$ a non-negative potential on
$\Sigma$, $K$ the Gaussian curvature and $a$ is a non-negative constant.

Such operators $L$ arise as the stability operator of $\Sigma$ immersed in a
Riemannian $3-$manifold with constant mean curvature (for particular choices of $V$
and $a$). We assume $L$ is non-positive acting on functions compactly supported on
$\Sigma$ and we obtain results in the spirit of some theorems of
Ficher-Colbrie-Schoen, Colding-Minicozzi and Castillon. We extend these theorems to
$a \leq 1/4$. We obtain results on the conformal type of $\Sigma$ and a distance (to
the boundary) lemma.
\end{abstract}

\section{Introduction}
A stable compact domain $\Sigma$ on a minimal surface in a Riemannian $3-$manifold
$M^3$, is one whose area can not be decreased up to second order by a variation of
the domain leaving the boundary fixed. Stable oriented domains $\Sigma$ are
characterized by the \emph{stability inequality} for normal variations $\psi N$
\cite{SY}

$$ \int _{\Sigma} \psi ^2 |A|^2 + \int _{\Sigma} \psi ^2 {\rm Ric}_{\mathcal{M}^3} (N,N)
\leq \int_{\Sigma} |\nabla \psi|^2  $$for all compactly supported functions $\psi
\in H^{1,2}_0 (\Sigma)$. Here $|A|^2$ denotes the the square of the length of the
second fundamental form of $\Sigma$, ${\rm Ric}_{\mathcal{M}^3} (N,N)$ is the Ricci
curvature of $\mathcal{M}^3$ in the direction of the normal $N$ to $\Sigma$ and
$\nabla $ is the gradient w.r.t. the induced metric.

One writes the stability inequality in the form
$$ \left.\frac{d^2}{dt^2}\right\vert _{t=0}{\rm Area}(\Sigma (t))=
- \int _{\Sigma} \psi L \psi \geq 0 ,$$where $L$ is the linearized operator of the
mean curvature
$$ L = \Delta + |A|^2 + {\rm Ric}_{\mathcal{M}^3} .$$

In terms of $L$, stability means that $-L$ is nonnegative, i.e., all its eigenvalues
are non-negative. $\Sigma$ is said to have finite index if $-L$ has only finitely
many negative eigenvalues.

Since the stability inequality is derived from the second variation formula under
normal variations of $\Sigma$, geometrically, $\Sigma$ has finite index if there is
only a finite dimensional space of normal variations which strictly decrease the
area.

In the $1970's$ and $80's$, this subject received an important number of
contributions (see \cite{B-doC,dCP,FCS,FC,Gu1,Gu2}), and even now it is a topic of
interest (see \cite{L,MPR} for surveys).

D. Fischer-Colbrie and R. Schoen \cite{FCS} studied stable surfaces by considering
the non-negativity of operators on a surface $\Sigma$ with a metric $g$ of the form
$$ L = \Delta + V - a K ,$$where $\Delta$ and $K$ are the Laplacian and Gaussian
curvature associated to $g$ respectively, $a$ is a positive constant and $V$ is a
nonnegative function. The main result of \cite{FCS} for stable surfaces in
three-manifolds is based on the following: \emph{For every complete metric on the
disk, there exists a number $a_0$ depending on the metric satisfying $0 \leq a_0 <
1$, so that for $a\leq a_0$, there is a positive solution of $\Delta -a K$, and for
$a>a_0$ there is no positive solution.} Note that the existence of a positive
solution of $\Delta -a K$ is equivalent to the non-positivity of $L =\Delta -a K$
(see \cite{FCS}).

Then a natural question was: \emph{What is the optimal $a_0$?} M. do Carmo and C. K.
Peng \cite{dCP} proved (implicitly) that $a_0 \geq 1/2$ for every complete metric on
the disk. Years later, S. Kawai \cite{K} (following ideas of A.V. Pogorelov
\cite{P}) stated that $a_0 \geq 1/4$ for a metric with non-positive curvature.

T. Colding and W. Minicozzi \cite{CM} introduced a new technique to study this type
of operator based on the first variation formula for length and the Gauss-Bonnet
formula. Using this technique they obtained a formula which, when $a>1/2$, gives
quadratic area growth of the geodesic disks on the surface and the integrability of
the potential $V$ at the same time (note that the stability operator can be realized
with the right choice of $V$).

Recently, P. Castillon \cite{Ca} used the ideas of Colding-Minicozzi to improve
their result to $a>1/4$. Moreover, he answered the question of the optimal value of
$a_0$ and proved the following:

{\bf Theorem:} \emph{Let $\Sigma$ be a complete noncompact Riemannian surface. Set
$$a_0 = {\rm sup}\set{a \in \r ^+ : \text{ there exists a positive solution to }
\Delta u - a K u =0 \text{ on } \Sigma} .$$ If $a_0 >1/4$ then $\Sigma$ is
conformally equivalent to $\c$ or $\c ^{*} = \c - \set{0}$.}

This improvement was by an appropriate choice of radial cut-off functions. In fact,
the same cut-off function has been considered in \cite{MPR} to obtain an analogous
formula to that of Colding-Minicozzi but for $a>1/4$.

Indeed, the value $a_0 =1/4$ is critical since this is the value for the
Poincar\'{e} metric on the disk. Thus one can not expect to have Castillon type
results without other hypothesis. That will be the main line of this work, the study
of these operators when $a\leq 1/4$ under other hypothesis on the area growth of the
geodesic disks. A subject  not studied until now, as far as we know.

The paper is organized as follows. In Section 2 we establish the notation that we
will use. We develop in Section 3 an inequality in the spirit of Colding-Minicozzi
for the quadratic form associated to the differential operator $L = \Delta + V - a K
$, $a>0$. We then apply this for a specific choice of radial cut-off function
defined on a geodesic disk on the surface, studying the behavior when the radius
tends to infinity. We will see why we need some hypothesis on the area growth of the
geodesic disks in the case that $0< a \leq 1/4$.

In Section 4 we consider the problem posed by D. Fischer-Colbrie and R. Schoen when
$0< a\leq 1/4$:

\vspace{.3cm}

{\bf Theorem \ref{t3}:}

{\it Let $\Sigma$ be a complete Riemannian surface with $k-$AAG. Set
$$a_0 = {\rm sup}\set{a \in \r ^+ : \text{ there exists a positive solution to }
\Delta u - a K u =0 \text{ on } \Sigma} .$$

If $k < 2$ and $0 \leq a_0 \leq 1/4$, then $\Sigma $ is conformally equivalent to
$\c$ or $\c ^{*}$. If $k=2$, $\Sigma$ is parabolic with finite topology.}

\vspace{.3cm}

Here, $k-$AAG means:

\vspace{.3cm}

{\bf Definition \ref{def1}:}

{\it Let $\Sigma$ be a Riemannian surface. We say that $\Sigma$ has Asymptotic Area
Growth of degree $k$ ($k-$AAG) if there exists positive constants $k ,C \in \r ^+ $
such that
$$ \lim _{r\To \infty}\frac{{\rm Area} ( D(p,r))}{r^k} = C  , \, \forall p \in \Sigma
.$$}

\vspace{.3cm}

In Section 5 we obtain a Huber-type Theorem. We prove

\vspace{.3cm}

{\bf Theorem \ref{t7}:}

{\it Let $\Sigma $ be a complete Riemannian surface with $k-$AAG and $0< a \leq
1/4$. Suppose $L = \Delta  - a K$ is non-positive on $\Sigma \setminus \Omega$,
$\Omega $ a compact set. Then, if $k \leq 2$, $\Sigma $ is conformally equivalent to
a compact Riemann surface with a finite number of points removed.}

\vspace{.3cm}

{\bf Theorem \ref{t8}:}

{\it Let $\Sigma $ be a complete noncompact parabolic Riemannian surface such that
$\int _{\Sigma}K^+ < +\infty$, with $K^{+}= {\rm max}\set{K,0}$. Suppose that $L =
\Delta - a K$ is non-positive on $\Sigma $, where $a> 0$. Then
\begin{itemize}
\item $K \in L^1 (\Sigma)$, i.e., it is integrable. In fact, $0 \leq \int _{\Sigma} K \leq 2 \pi \chi (\Sigma)
$.
\item $\Sigma $ has quadratic area growth.
\item $\Sigma$ is conformally equivalent either to the plane or to the cylinder.
\end{itemize}}

\vspace{.3cm}

We will apply these results on Section 6 to stable surfaces, we will be able to
bound the distance of any point to the boundary, this is known as the \emph{Distance
Lemma} (see \cite{R2}, \cite{RR} or \cite{MPR} in the more general version, i.e. for
$a>1/4$ and $V\equiv c >0$ some constant, and \cite{Ma} for a sharp bound in space
forms).

In fact, the authors gave an explicit bound of this distance in terms of $a$, when
$a>1/4$, and $c>0$. Here we generalize this result for $0<a \leq 1/4$, giving the
existence of some constant which bounds this distance.

\vspace{.3cm}

{\bf Theorem \ref{t5}:}

{\it Let $\Sigma $ be a Riemannian surface possibly with boundary and $k-$AAG.
Suppose that $L = \Delta + V - a K$ is non-positive, where $V \geq c >0 $ and $0< a
\leq 1/4$. Then, there exists a positive constant $C$ such that
$$ {\rm dist}_{\Sigma}(p , \partial \Sigma) \leq C , \, \, \forall p \in \Sigma .$$

In particular, if $\Sigma$ is complete with $\partial \Sigma = \emptyset$ then it
must be topologically a sphere.}

\vspace{.3cm}

In addition, we will control the growth of the integral of the potential (known for
$a>1/4$); that is,

\vspace{.3cm}

{\bf Theorem \ref{t6}:}

{\it Let $\Sigma $ be a complete Riemannian surface satisfying $k-$AAG, $k\leq 2$.
Suppose that $L = \Delta + V - a K$ is non-positive, where $V \geq 0 $ and $0< a
\leq 1/4$. Then, $V \in L^1(\Sigma)$, i.e., $V$ is integrable.

Moreover, if $\Sigma$ has $k-$AAG with $k>2$, then for $2(b+1) \geq  k$ we have
\begin{equation*}
\int _{D(s)} V \leq C s^{2b}
\end{equation*}for some positive constant $C$.}

Finally, we consider a problem posed in \cite{FCS} for stable surfaces immersed in a
three-manifold. In \cite[Theorem 3]{FCS}, they proved: \emph{Let $N$ be a complete
oriented $3-$manifold of non-negative scalar curvature. Let $\Sigma$ be an oriented
complete stable minimal surface in $N$. If $\Sigma$ is noncompact, conformally
equivalent to the cylinder and the absolute total curvature of $\Sigma$ is finite,
then $\Sigma$ is flat and totally geodesic.}

And they state \cite[Remark 2]{FCS}: \emph{We feel that the assumption of finite
total curvature should not be essential in proving that the cylinder is flat and
totally geodesic.}

Using Theorem \ref{t8}, we are able to partially answer this question.

\vspace{.3cm}

{\bf Theorem \ref{t9}:}

{\it Let $N$ be a complete oriented $3-$manifold of non-negative scalar curvature.
Let $\Sigma$ be an oriented complete stable minimal surface in $N$. If $\Sigma$ is
noncompact, conformally equivalent to the cylinder and $\int _{\Sigma} K^{+}$ is
finite, then $\Sigma$ is flat and totally geodesic.}

\section{Preliminaries}

We denote by $\Sigma$ a connected Riemannian surface, with riemannian metric $g$,
and possibly with boundary $\partial \Sigma$. Let $p_0 \in \Sigma$ be a point of the
surface and $D(p_0,s)$, for $s>0$, denote the geodesic disk centered at $p_0$ of
radius $s$. We assume that $\overline{D(p_0 ,s)} \cap \partial \Sigma = \emptyset$.
Moreover, let $r$ be the radial distance of a point $p$ in $D(p_0, s)$ to $p_0$. We
write $D(s)=D(p_0 ,s)$.

We also denote
\begin{eqnarray*}
l(s) &=& {\rm Length}(\partial D(s)) \\
a(s) &=& {\rm Area}(D(s))\\
K(s) &=& \int _{D(s)} K \\
\chi (s)&=& \text{Euler characteristic of } D(s)
\end{eqnarray*}

Moreover, we will need the following result due to K. Shiohama and M. Tanaka (see
\cite{ST}) which follows from the first variation formula for length and the
Gauss-Bonnet formula,

\begin{theorem}\label{t1}
The function $l$ is differentiable almost everywhere and we have
\begin{enumerate}
\item for almost all $r \in \r$,
\begin{equation}\label{derlength}
l'(r) \leq 2 \pi \chi (r)- K(r),
\end{equation}
\item for all $0 \leq a < b$,
\begin{equation}\label{intlength}
l(b)-l(a) \leq \int_ a ^b l' (r)
\end{equation}
\end{enumerate}
\end{theorem}

Here, $'$ denotes the derivative with respect to $r$.

Let $L = \Delta + V - a K$ be a differential operator on $\Sigma$ acting on
compactly supported $f \in H^{1,2}_0 (\Sigma)$, where $a >0$ is constant, $V \geq
0$, $\Delta$ and $K$ the Laplacian and Gauss curvature associated to the metric $g$
respectively.

The index form of these kind of operators is
\begin{equation}\label{varcharac}
I(f)=\int _{\Sigma } \left\{ \norm{\nabla f}^2 - V f^2+ a K f^2 \right\}
\end{equation}where $\nabla $ and $\| \cdot \|$ are the gradient and norm associated
to the metric $g$. One has
$$\int_{\Sigma} f L f = -I(f) .$$

We will use the following condition on the area growth of $\Sigma$

\begin{definition}\label{def1} Let $\Sigma$ be a Riemannian surface. We
say that $\Sigma$ has Asymptotic Area Growth of degree $k$ ($k-$AAG) if there exists
positive constants $k ,C \in \r ^+ $ such that
$$ \lim _{r\To \infty}\frac{{\rm Area} ( D(p,r))}{r^k} = C  , \, \forall p \in \Sigma .$$
\end{definition}

Note that, by the Triangle Inequality, this condition does not depend on the point
$p$.

\section{A Colding-Minicozzi stability inequality}

Here, we will establish a general inequality for $I(f)$ (see \eqref{varcharac}) when
$f $ is a radial function defined on a geodesic disk, following the method used by
T. Colding and W. Minicozzi in \cite{CM}. The proof of this can be found in
\cite{Ca}, but we include it here for the sake of completeness. The final
formulation is slightly different than that of \cite{Ca}.

\begin{lemma}[Colding-Minicozzi stability inequality]\label{l1}
Let $\Sigma $ be a Riemannian surface possibly with boundary and $K \not \equiv 0$.
Let us fix a point $p_0\in \Sigma$ and positive numbers $0 \leq \eps < s $ such that
$\overline{D(s)}\cap \partial \Sigma = \emptyset$. Let us consider the differential
operator $L = \Delta + V - a K$, where $V \geq 0$ and $a$ is a positive constant,
acting on $f \in H^{1,2}_0 (\Sigma)$. Let $f: D(s) \To \r $ a non-negative radial
function, i.e. $f \equiv f(r)$, such that

$$\begin{matrix}
f(r) \equiv  1, & \text{ for } r \leq \eps \\
f(r) \equiv  0, & \text{ for } r \geq s \\
f'(r) \leq 0 , & \text{ for } \eps < r < s
\end{matrix}$$

Then, the following holds
\begin{equation}\label{CMineq}
 I(f) \leq 2a \left( \pi G (s) - f'_{-}(\eps) l(\eps)\right) -\int _{D(s)} V f(r)^2
 + \int _{\eps}^s \left\{ (1-2a)f'(r) ^2 - 2 a f(r)f''(r)\right\}l(r) ,
\end{equation}where
\begin{eqnarray*}
G(s) &:= &- \int _0 ^s (f(r)^2)' \chi (r) \leq 1 ,\\
f'_{-}(\eps) &:=& \lim _{\begin{matrix} r \to \eps \\ \eps < r\end{matrix}} f'(r) .
\end{eqnarray*}
\end{lemma}
\begin{proof}
Let us denote
$$ \alpha = \int _{D(s)} \norm{\nabla f}^2 , \, \, \beta = \int _{D(s)} K f^2 .$$

On the one hand, by the Co-Area Formula
\begin{equation*}
\alpha = \int _{D(s)} \norm{\nabla f}^2 = \int _{\eps}^s f'(r)^2 \int _{\partial
D(r)}1 = \int _{\eps}^s f'(r)^2 l(r)
\end{equation*}

On the other hand, by Fubini's Theorem and integrating by parts, we have
\begin{equation*}
\beta = \int _0 ^s f(r)^2 \int _{\partial D(r)} K = \int _0 ^s f(r)^2 K'(r) = - \int
_0 ^s (f(r)^2)' K(r) .
\end{equation*}

Now, by \eqref{derlength} and $(f(r)^2)' = 2 f(r)f'(r)\leq 0$, we have
$$ -(f(r)^2)' K(r) \leq (f(r)^2)' (l'(r) - 2 \pi \chi (r)) .$$

Integrating by parts and taking into account that $\int _0 ^s (f(r)^2)' = -1$, we
obtain
\begin{equation*}
\begin{split}
\beta &  \leq \int _0 ^s (f(r)^2)' (l'(r) -2\pi \chi (r)) = -2\pi  \int _0 ^s
(f(r)^2)' \chi (r) + \int _0 ^s (f(r)^2)' l'(r)\\
 & = 2\pi G(s)+ \int  _0 ^s (f(r)^2)' l'(r) = 2\pi G(s)+ \int_{\eps} ^s ((f(r)^2)' l(r))'
 - \int _{\eps} ^s (f(r)^2)'' l(r) \\
 &= 2  \pi G(s)  - 2 f'_{-}(\eps)l(\eps) - \int _{\eps} ^s (f(r)^2)'' l(r) .
\end{split}
\end{equation*}

Thus, putting $\alpha$ and $\beta$ together
\begin{equation*}
\int _{D(s)} \left\{ \norm{\nabla f}^2 + a K f^2 \right\} \leq 2a \left( \pi G(s) -
f'_{-}(\eps) l(\eps)\right) + \int _{\eps}^s \left\{ (1-2a)f'(r) ^2 - 2a
f(r)f''(r)\right\}l(r) .
\end{equation*}

Note that the bound on $G(s)$ follows since the Euler characteristic of $D(s)$ is
less than or equal to $1$.
\end{proof}

Now, we will work with the special radial function given by

\begin{equation}\label{function}
f(r)=\left\{\begin{matrix}
1 & r \leq s e^{-s}  \\[3mm]
\left( \dfrac{\ln (s/r)}{s}\right) ^ b & s e^{-s} \leq r \leq s \\[3mm]
0 & r \geq s
\end{matrix}\right.
\end{equation}where $s>0$, $b\geq 1$ and $r$ is the radial distance of a point $p$ in $D(s)$ to $p_0$.

We summarize the properties of this function in the following result

\begin{proposition}\label{p1}
Let $a$ and $s$ be positive constants and $f: [0,s] \To \r $ the function given by
\eqref{function}. Denote
\begin{eqnarray}
\alpha &=&1+ b \frac{1-4a}{2a}\label{alpha}\\[3mm]
g(r) &=& \frac{\ln (s/r)}{s} \label{g}\\[3mm]
\phi (r) &=& \alpha - s g(r) \label{phi} \\[3mm]
F(r) &=& (1-2a)f'(r) ^2 - 2a f(r)f''(r) \label{F}
\end{eqnarray}

Then, for $r \in (s e^{-s} , s)$, we have
\begin{eqnarray}
f'(r) &=& - \frac{b}{s r}g(r)^{b-1} \leq 0   \label{derf}\\[3mm]
F(r) &=& 2 ab \frac{g(r)^{2(b-1)}}{s^2 r^2} \phi (r) \label{Fr}
\end{eqnarray}

Moreover, if $\alpha > 0$ and $s  > \alpha + \delta > \alpha >0$ for some positive
constant $\delta$, then, the intervals
$$ \mathcal{I}_1 = [s e^{-s}, s e^{-(\alpha +\delta)}], \, \mathcal{I}_2 = [s e^{-(\alpha +\delta)}, s e^{-\alpha}], \,
\mathcal{I}_3 = [s e^{-\alpha}, s ] $$are well defined and
\begin{eqnarray}
F _{\vert \mathcal{I}_1} &\leq & -2\delta ab (\alpha + \delta)^{2(b-1)}e^{2(\alpha +\delta)}\frac{1}{s^{2(b+1)}} \label{Fint1pos}\\[3mm]
F _{\vert \mathcal{I}_2} &\leq & 0 \label{Fint2pos}\\[3mm]
F _{\vert \mathcal{I}_3} &\leq & 2 ab \alpha ^{2b -1}e^{2\alpha}\frac{1}{s^{2(b+1)}}
\label{Fint3pos}
\end{eqnarray}

\end{proposition}
\begin{proof}
First, \eqref{derf} and \eqref{Fr} are straightforward computations using the
definitions of \eqref{alpha}, \eqref{g}, \eqref{phi} and \eqref{F}.

Let us assume that $\alpha >0$. Let $s>0$ be a positive number such that $s > \alpha
+ \delta$ for some $\delta >0 $ fixed. Then
$$ e^{-s} < e^{-(\alpha + \delta)} < e^{-\alpha} < 1 $$which means that the
intervals $\mathcal{I}_i$, $i=1,2,3$, are well defined.

Since $g$ (given by \eqref{g}) is a decreasing function, we have
\begin{equation*}
\begin{matrix}
\dfrac{\alpha + \delta}{s } & \leq &  g_{\vert \mathcal{I}_1} & \leq & 1 \\[3mm]
\dfrac{\alpha }{s } & \leq &  g_{\vert \mathcal{I}_2} & \leq & \dfrac{\alpha +
\delta}{s }
\\[3mm]
0 & \leq &  g_{\vert \mathcal{I}_3} & \leq & \dfrac{\alpha }{s }
\end{matrix} .
\end{equation*}

Thus,
\begin{equation*}
\begin{matrix}
  \dfrac{(\alpha + \delta)^{2 (b-1)} e^{2(\alpha + \delta)}}{s^{2b} } & \leq &
 \left(\dfrac{g(r)^{2(b-1)}}{r^2} \right)_{\vert \mathcal{I}_1} & \leq & \dfrac{e^{2s}}{s^2} \\[6mm]
 \dfrac{\alpha ^{2(b-1)} e^{2\alpha}}{s^{2b} } & \leq &
\left(\dfrac{g(r)^{2b)}}{r^2} \right)_{\vert \mathcal{I}_2} & \leq &\dfrac{(\alpha +
\delta)^{2(b-1)} e^{2(\alpha + \delta)}}{s^{2b} }
\\[6mm]
0 & \leq &  \left(\dfrac{g(r)^{2(b-1)}}{r^2} \right)_{\vert \mathcal{I}_3} & \leq &
 \dfrac{\alpha ^{2(b-1)} e^{2\alpha}}{s^{2b} }
\end{matrix}
\end{equation*}and $\phi$ (given by \eqref{phi}) satisfies
\begin{equation*}
\begin{matrix}
\alpha - s & \leq &  \phi _{\vert \mathcal{I}_1} & \leq & -\delta \\[3mm]
-\delta & \leq &  \phi _{\vert \mathcal{I}_2} & \leq & 0 \\[3mm]
0 & \leq &  \phi _{\vert \mathcal{I}_3} & \leq & \alpha
\end{matrix}
\end{equation*}

Hence, from \eqref{Fr},
\begin{eqnarray*}
F _{\vert \mathcal{I}_1} &\leq & -2\delta ab (\alpha + \delta)^{2(b-1)}e^{2(\alpha +\delta)}\frac{1}{s^{2(b+1)}} \\[3mm]
F _{\vert \mathcal{I}_2} &\leq & 0 \\[3mm]
F _{\vert \mathcal{I}_3} &\leq & 2 ab \alpha ^{2b -1}e^{2\alpha}\frac{1}{s^{2(b+1)}}
;
\end{eqnarray*}as desired.
\end{proof}

Given $0 \leq r_1 < r_2 $, we denote
$$ A(r_1 ,r_2) = a(r_2)-a(r_1) ,$$that is,
$$ A(r_1, r_2) = \int _{r_1}^{r_2} l(r) .$$

Also, we denote
$$ \mathcal{K}(r_1) = {\rm min}_{[0,r_1]}\set{K(r)} . $$

\begin{lemma}\label{l2}
Let $\Sigma $ be a Riemannian surface possibly with boundary and $K \not \equiv 0$.
Let us fix a point $p_0\in \Sigma$ and a positive number $ s > 0$ such that
$\overline{D(s)}\cap \partial \Sigma = \emptyset$. Set $L = \Delta + V - a K$, where
$V \geq 0$ and $a$ is a positive constant, acting on $f \in H^{1,2}_0 (\Sigma)$.
Given $b\geq 1$, let $\alpha $ be defined by \eqref{alpha}. Then, if $\alpha > 0$,
\begin{equation}\label{estimatepos}
 I(f) \leq 2a \left(G(s) \pi + b\frac{2\pi - \mathcal{K}(se^{-s})}{s}\right)+ \rho^{+}_{a,b}(\delta ,s) -
\left(\int _{D(s e^{-s})}  V  + \int _{D(s)\setminus D(s e^{-s})} \left(
\frac{\ln(s/r)}{s}\right) ^{2b} V \right)
\end{equation}where
\begin{equation}\label{rhopos}
\rho ^{+} _{a,b}(\delta ,s) =2ab \alpha ^{2b-1}e^{2\alpha}\left(
\frac{A(se^{-\alpha} ,s)}{s^{2(b+1)}} - \frac{\delta e^{2\delta}}{\alpha}\left(
1+\frac{\delta}{\alpha}\right)^{2(b-1)}\frac{A(se^{-s},
 se^{-(\alpha+\delta)})}{s^{2(b+1)}}\right)
\end{equation}
\end{lemma}

\begin{proof}
We will use the function $f$ given by \eqref{function} in the equation
\eqref{CMineq} of Lemma \ref{l1} taking into account that $s e^{-s}$ plays the role
of $\eps$ in the formula, i.e. $\eps = s e^{-s}$.

First, we will estimate the term $f'(\eps) l (\eps)$ in \eqref{CMineq}.

Using \eqref{derlength} and \eqref{intlength}, for any $\eps >0$ we have
\begin{equation}\label{leps}
l(\eps) \leq \int _0 ^{\eps} l'(r) \leq 2\pi \eps - {\rm min}_{r\in
[0,\eps]}\set{K(r)} \eps = \left( 2 \pi - \mathcal{K}(\eps)\right) \eps .
\end{equation}

Also, by \eqref{derf},
$$ f'(s e^{-s}) = -\frac{b}{s^2 e^{-s}} ,$$so, from \eqref{leps} we obtain
$$ - f'(s e^{-s}) l(s e^{-s}) = b \frac{1}{s}\frac{l(s e^{-s})}{s e^{-s}} \leq
b \frac{2 \pi - \mathcal{K}(s e^{-s})}{s}.$$

Thus,
\begin{equation}\label{term1}
2a\left( G(s) \pi - f'(s e^{-s})l(s e^{-s})\right) \leq 2a \left(G(s) \pi +
b\frac{2\pi - \mathcal{K}(se^{-s})}{s}\right) .
\end{equation}

Now, note that with the notation of Proposition \ref{p1}, we have

$$\int _{s e^{-s}} ^{s} \left( (1-2a)f'(r) ^2 - 2a f(r)f''(r)\right) l(r) =
\int _{s e^{-s}}^{s} F(r)l(r).$$

Thus, from \eqref{Fint1pos}, \eqref{Fint2pos} and \eqref{Fint3pos}
\begin{equation*}
\begin{split}
\int _{s e^{-s}}^{s} F(r) l(r) &\leq \int _{\mathcal{I}_1} F(r)l(r) + \int
_{\mathcal{I}_2} F(r)l(r) +\int _{\mathcal{I}_3} F(r)l(r) \\
 & \leq \int _{\mathcal{I}_1} F(r)l(r) + \int _{\mathcal{I}_3} F(r)l(r) \\
 &\leq 2 ab \alpha ^{2b-1}e^{2\alpha} \frac{1}{s^{2(b+1)}}\int _{\mathcal{I}_3}l(r)
  - 2 ab \delta \left( \alpha + \delta\right)^{2(b-1)}e^{2(\alpha+\delta)}
  \frac{1}{s^{2(b+1)}} \int _{\mathcal{I}_1} l(r)\\
  & =  2 ab \alpha ^{2b-1}e^{2\alpha} \left( \frac{A(se^{-\alpha} ,s)}{s^{2(b+1)}} -
 \frac{\delta e^{2\delta}}{\alpha}\left( 1+\frac{\delta}{\alpha}\right)^{2(b-1)}\frac{A(se^{-s},
 se^{-(\alpha+\delta)})}{s^{2(b+1)}}\right),
\end{split}
\end{equation*}that is
\begin{equation}\label{term2pos}
\int _{s e^{-s}}^{s} F(r) l(r) \leq 2 ab \alpha ^{2b-1}e^{2\alpha} \left(
\frac{A(se^{-\alpha} ,s)}{s^{2(b+1)}} -
 \frac{\delta e^{2\delta}}{\alpha}\left( 1+\frac{\delta}{\alpha}\right)^{2(b-1)}\frac{A(se^{-s},
 se^{-(\alpha+\delta)})}{s^{2(b+1)}}\right),
\end{equation}hence, combining \eqref{term1} and \eqref{term2pos}, we obtain
\eqref{estimatepos}.
\end{proof}

\begin{remark}
Thus, it is clear from the above Lemma that the behavior of $I(f)$ depends on the
function $\rho ^{+}$, which depends on the area growth of the surface.
\end{remark}

The last results in this Section are devoted to the asymptotic behavior of the
function $\rho ^{+}$ under suitable conditions on the surface.

\begin{lemma}\label{l3}
Let $\Sigma $ be a Riemannian surface possibly with boundary satisfying $k-$AAG and
$K \not \equiv 0$. Given $b\geq 1$ and $a>0$, let $\alpha $ be defined by
\eqref{alpha}. Then, if $\alpha >0$, the asymptotic behavior of $\rho ^{+}_{a,b}$,
given by \eqref{rhopos}, as $s $ goes to infinity is
\begin{equation}\label{asymppos}
\rho ^{+} _{a,b} (\delta ,s ) \sim C^{+} \frac{s^k}{s^{2(b+1)}}
\tilde{\rho}^{+}_{\alpha ,k}(\delta)
\end{equation}where
\begin{eqnarray}
C^{+}(a,b,C) &=& 2 a b C \alpha ^{2b -1} e^{2\alpha} \label{Cpos}\\
\tilde{\rho}^{+}_{\alpha ,k}(\delta) &=& 1 - e^{-k \alpha}\left( 1+ e^{(2-k)\delta}
\frac{\delta}{\alpha}\left(1+\frac{\delta}{\alpha}\right)^{2(b-1)}\right)
\label{rtildepos}
\end{eqnarray}and $C$ is the positive constant such that $a(s) \sim C s^k$.
\end{lemma}
\begin{proof}
We want to control the asymptotic behavior of $\rho ^+ $. Since $\Sigma$ has
$k-$AAG, this means that there exists $C >0$ such that
$$ a(s) \sim C s^k $$for $s$ large.

Hence
\begin{eqnarray*}
A(s e^{-s},s e^{-(\alpha + \delta)}) &\sim & C s^k e^{-k(\alpha + \delta)} \\
A(s e^{-\alpha} ,s) & \sim & C s^k (1- e^{-k\alpha})
\end{eqnarray*}thus, from \eqref{rhopos},

\begin{equation*}
\begin{split}
\rho ^+_{a,b}(s ,\delta) & \sim 2 ab C \alpha ^{2b -1} e^{2\alpha} \frac{s
^k}{s^{2(b+1)}} \left( 1- e^{-k\alpha}  - e^{2\delta}e^{-k(\alpha
+\delta)}\frac{\delta}{\alpha}\left( 1+\frac{\delta}{\alpha}\right)^{2(b-1)}\right)
\\
 & = \left( 2 a b C \alpha ^{2b -1} e^{2\alpha} \right)\left( 1 - e^{-k \alpha}\left( 1+ e^{(2-k)\delta}
\frac{\delta}{\alpha}\left(1+\frac{\delta}{\alpha}\right)^{2(b-1)}\right)\right)
\frac{s ^k}{s^{2(b+1)}} \\
 & = C^+ \frac{s ^k}{s^{2(b+1)}} \tilde{\rho}^{+}_{\alpha ,k}(\delta)
\end{split}
\end{equation*}as desired.
\end{proof}

\begin{remark}\label{rdelta0}
Let us note that the behavior of $\tilde{\rho}^+$ depends strongly on the degree of
the AAG. Moreover, we have that $\tilde{\rho}^{+}_{\alpha , k} (\delta)$ is a
bounded function of $\delta \in \r ^+ $ since it is continuous and
$$ \lim _{\delta \To 0} \tilde{\rho}^{+}_{\alpha , k}
(\delta) = 1- e^{-k\alpha }=  \lim _{\delta \To +\infty} \tilde{\rho}^{+}_{\alpha ,
k} (\delta).$$

Thus there exists $\delta _0 >0 $ such that
$$\rho _{{\rm min}}= \tilde{\rho} _{\alpha , k}^+ (\delta _0) =
{\rm min}_{\delta >0} \set{\tilde{\rho}^+_{\alpha ,k}}.$$
\end{remark}

So, with this last remark in mind, we conclude

\begin{corollary}\label{c1}
Assuming the conditions of Lemma \ref{l3}, if $\alpha >0$ and $2(b+1)>k$, then as $s
\to + \infty$
\begin{eqnarray}\label{asymk2pos}
\rho^{+}_{a ,b}(\delta _0 , s) \to 0
\end{eqnarray}where $\delta _0$ is given in Remark \ref{rdelta0}.
\end{corollary}

\section{On a problem of D. Fischer-Colbrie and R. Schoen}

In \cite{FCS}, the authors proved: \emph{For every complete metric on the disc,
there exist a number $a_0$ depending on the metric satisfying $0 \leq a_0 < 1$ so
that for $a\leq a_0$ there is a positive solution of $\Delta -a K$, and for $a>a_0$
there is no positive solution.} Here, $\Delta$ and $K$ denote the laplacian and
Gauss curvature of the metric respectively.

As we said in the Introduction, P. Castillon \cite{Ca} proved the following:

{\bf Theorem:} \emph{Let $\Sigma$ be a complete noncompact Riemannian surface. Set
$$a_0 = {\rm sup}\set{a \in \r ^+ : \text{ there exists a positive solution to }
\Delta u - a K u =0 \text{ on } \Sigma} .$$ If $a_0 >1/4$ then $\Sigma$ is
conformally equivalent to $\c$ or $\c ^{*} = \c - \set{0}$.}

The method used for this is a formula as in Lemma \ref{l1} (to control the conformal
type of the ends).

Moreover, in \cite{Ca} and \cite{MPR}, it is shown that if $L_a = \Delta - a K \leq
0$ and $a>1/4$, then $\Sigma$ has at most quadratic area growth, i.e.,
$$ a(s) \leq C s^2 $$for some positive constant $C$ and all $s>0$.

But, assuming some $k-$AAG on the surface we obtain the following (this is the first
result we know of when $a_0 \in [0,1/4]$).

\begin{theorem}\label{t3}
Let $\Sigma$ be a complete Riemannian surface with $k-$AAG. Set
$$a_0 = {\rm sup}\set{a \in \r ^+ : \text{ there exists a positive solution to }
\Delta u - a K u =0 \text{ on } \Sigma} .$$

If $k < 2$ and $0 \leq a_0 \leq 1/4$, then $\Sigma $ is conformally equivalent to
$\c$ or $\c ^{*}$. If $k=2$, $\Sigma$ is parabolic with finite topology.
\end{theorem}
\begin{proof}
Suppose that there exists $0<a \leq 1/4$ such that there exists a positive $u$
solution to $\Delta u - a K u = 0$ on $\Sigma$, then $L = \Delta - a K$ is non
positive.

On the one hand, consider the radial function $f(r)= (1-r/s)$, $r \leq s$, in the
equation \eqref{CMineq} of Lemma \ref{l1} with $V \equiv 0$, then we obtain:

$$0 \leq I(f) \leq 2a \pi G(s) + \frac{(1-2a)}{s^2}a(s)$$

On the other hand, assume that there exists $s_0$ so that for $s \geq s_0$ we have
$\chi (s)\leq -M $, and then
\begin{equation*}
\begin{split}
G(s) & = - \int _0 ^s (f(r)^2)' \chi (r) = - \int _0 ^{s_0} (f(r)^2)' \chi (r)- \int
_{s_0} ^s (f(r)^2)' \chi (r)\\
 & \leq - \int _0 ^{s_0} (f(r)^2)' + M \int _{s_0} ^s (f(r)^2)' =  -\left( f(s_0) ^2 - f(0)^2\right) +
 M \left( f(s)^2 - f(s_0)^2 \right) \\
 &= - (M+1) f(s_0)^2 +1 = -(M+1)\left( 1-s_0/s\right)^2 +1 ,
\end{split}
\end{equation*}

Thus, joining the above inequalities, we obtain
$$ 0 \leq 2 a \left( 1- (M+1)(1-s_0/s)^2\right) + \frac{1-2a}{s^2}a(s). $$

Set $$ T(s) := 2 a \left( 1- (M+1)(1-s_0/s)^2\right) + \frac{1-2a}{s^2}a(s) ,$$so
$$  0 \leq \lim_{s \to + \infty} T(s) = - M + (1-2a)\lim_{s \to +\infty}\frac{a(s)}{s^2} .$$

If $k< 2$, then $M \geq 0$, this means that $\Sigma $ is homeomorphic either to a
plane or a cylinder.

If $k =2$, we have
$$0 \leq \lim _{s \to + \infty} T(s) = - M + (1-2a)C ,$$for some positive constant
$C$. Then, $M \leq (1-2a)C$, that is, $\Sigma $ has finite topology.

It remains to prove that each end of $\Sigma$ is parabolic. This is clear since the
area growth is quadratic.
\end{proof}

\begin{remark}
Even with the AAG hypothesis, this result is sharp. R. Schoen pointed out to us that
there exist hyperbolic surfaces with polynomial area growth bigger that 2. Let us
explain this. Consider the rotationally symmetric metric
$$ g = dr ^2 + \tau (r)^2 d\theta $$on $\r ^2$, with $\tau (r)= \frac{r^{1+\eps}}{\eps}$ for $\eps
>0$. Then, it is easy to see that the area of the geodesic disks are given by
$$ {\rm Area}(D(r)) = \frac{2\pi}{\eps (1+\eps)}r^{2+\eps} ,$$and the Gaussian
curvature is
$$ K = -\frac{\tau ''}{\tau } = - \frac{1+\eps}{r^2} .$$

Now, for $r$ large, we have
$$K \leq - \frac{1+\eps}{r^2 \ln r}$$and $\tau $ is unbounded. Hence, using \cite[Theorem
1]{M}, $(\r ^2 , g)$ is conformally hyperbolic.
\end{remark}

\section{A Huber-type Theorem and parabolicity}

Here we will establish a Huber type Theorem for surfaces with $2-$AAG and $L =
\Delta  - a K$, $0< a \leq 1/4$. In fact, the proof follows from the work of P.
Castillon \cite{Ca}.

\begin{theorem}\label{t7}
Let $\Sigma $ be a complete Riemannian surface with $k-$AAG and $0< a \leq 1/4$.
Suppose $L = \Delta  - a K$ is non-positive on $\Sigma \setminus \Omega$, $\Omega $
a compact set. Then, if $k \leq 2$, $\Sigma $ is conformally equivalent to a compact
Riemann surface with a finite number of points removed.
\end{theorem}
\begin{proof}
The main steps in \cite[Theorem B]{Ca} are controlling the topology and the area
growth of the surface. Note that once we know that the surface has at most quadratic
area growth, we control the conformal type of the ends. So, in the first item, as we
are assuming at most $2-$AAG, this last part is guaranteed. So, it remains to prove
that the topology is finite.

We follow the proof \cite[Proposition 3.1]{Ca}.

Let $s_0 ,s_1 > 0 $ such that $\Omega \subset D(s_0)$ and $s_0 < s_1 -1$. Define
$f_0 : [s_1 -1 , s_1] \To \r$ by $f_0(r) = r-s_1 +1$ and
$$ c_{a} =  - a K(s_1) +\int _{D(s_1)\setminus D(s_1 -1)} \left\{ \norm{\nabla f_0 (r)}^2 + a K f_0(r)^2\right\} ,
$$which is a constant depending on $a$ and the metric.

Now, consider the radial function
\begin{equation*}
f(r)=\left\{\begin{matrix}
f_0 (r) & \text{for} & r \in [s_1 -1 , s_1] \\[3mm]
1 & \text{for} & r \in [s_1 , s_2 ] \\[3mm]
\dfrac{s-r}{s-s_2}& \text{for} & r \in [s_2 , s]\\[3mm]
0 & \mbox{} & \text{elsewhere}
\end{matrix}\right.
\end{equation*}

Note that $f$ has compact support on $\Sigma \setminus D(s_0)$, so applying that $L$
is non positive, and following the computations of Lemma \ref{l1}, we obtain
\begin{equation}\label{huber}
0 \leq c_a + 2\pi a G(s) + \frac{2a l(s_2)}{s-s_2} + \frac{1-2a}{(s-s_2)^2}\int
_{s_2}^s l(r) ,
\end{equation}where
\begin{equation*}
\begin{split}
G(s) &= - \int _{s_2} ^s \left( \left(\frac{s-r}{s-s_2}\right)^2 \right)' \chi (r) .
\end{split}
\end{equation*}

Since $\Sigma$ has at most $2-$AAG, we have that
\begin{equation*}
c_a +\frac{2a l(s_2)}{s-s_2} + \frac{1-2a}{(s-s_2)^2}\int _{s_2}^s l(r) \To C
\end{equation*}as $s \To +\infty$, $C$ a positive constant.

If $\Sigma$ has infinite topology, then
$$ \liminf _{s \To + \infty} \chi (s) = - \infty ,$$that is, we can take $s_2$ big
enough so that $\chi (s) \leq -\frac{C+1}{2a\pi}$ for all $s \geq s_2$, therefore
\begin{equation*}
G(s) \leq \frac{C+1}{2a \pi} \int _{s_2}^s \left( \left(\frac{s-r}{s-s_2}\right)^2
\right)' \leq - \frac{C+1}{2a\pi} .
\end{equation*}

And so
$$0 \leq 0 \leq c_a + 2\pi a G(s) + \frac{2a l(s_2)}{s-s_2} +
\frac{1-2a}{(s-s_2)^2}\int _{s_2}^s l(r) \to -1 $$as $s \to + \infty$, which is a
contradiction. This completes the proof.
\end{proof}

In the first two sections we obtained parabolicity from the area growth of the
surface, but it is interesting (as we will see in the next Section ) to study what
happens when we assume parabolicity but not $k-$AAG.

\begin{theorem}\label{t8}
Let $\Sigma $ be a complete noncompact parabolic Riemannian surface such that $\int
_{\Sigma}K^+ < +\infty$, with $K^{+}= {\rm max}\set{K,0}$. Suppose that $L = \Delta
- a K$ is non-positive on $\Sigma $, where $a> 0$. Then
\begin{itemize}
\item $K \in L^1 (\Sigma)$, i.e., it is integrable. In fact, $0 \leq \int _{\Sigma} K \leq 2 \pi \chi (\Sigma)
$.
\item $\Sigma $ has quadratic area growth.
\item $\Sigma$ is conformally equivalent either to the plane or to the cylinder.
\end{itemize}
\end{theorem}
\begin{proof}
The two last statements follow from the first one; let us explain this briefly.

Assume that $K \in L^1 (\Sigma)$, then \eqref{derlength} implies that
$$ l'(r) \leq 2\pi - K(r) \leq C ,$$for some positive constant $C$, which means that
$l(r) \leq C r$ from \eqref{intlength}. Thus, $\Sigma$ has at most quadratic area
growth. Now, either Theorem \ref{t3} for $a<1/4$ or \cite[Theorem A]{Ca} gives us
the conformal type of the surface.

For a fixed point $p \in \Sigma$ and a sequence $s_0 < s_1 < s_2 < \ldots  \to +
\infty$, let us consider the sequence of positive functions defined by
\begin{equation*}
\begin{matrix}
\Delta f_i = 0 & \text{ on } & D(s_i)\setminus \overline{D(s_0)}\\
f_i = 1 & \text{ on } & \overline{D (s_0)}\\
f_i = 0 & \text{ on } & \partial D (s_i)
\end{matrix}.
\end{equation*}

Moreover, this sequence converges uniformly on compact subsets of $\Sigma $ to the
constant function $1$, and also is a monotone sequence by the Maximum Principle (see
\cite[Lemma 3.6]{L}).

So, following \cite[Theorem 10.1]{L}, using the boundary conditions and the fact
that $f_i$ is harmonic on $D(s_1)\setminus D(s_0)$, we have
\begin{equation*}
\begin{split}
\int _{D(s_i)\setminus D(s_0)} \norm{\nabla f_i}^2 &= \int _{\partial D(s_i)} f_i
\frac{\partial f_i}{\partial \eta} - \int _{\partial D(s_0)} f_i \frac{\partial
f_i}{\partial \eta} \\
 & = - \int _{\partial D(s_0)} \frac{\partial f_i}{\partial \eta}
\end{split}
\end{equation*}where $\frac{\partial f_i}{\partial \eta}$ the outward pointing
derivative.

Thus, using that $f_ i \To 1$ (uniformly on compact subsets), the right hand must
goes to $0$ as $s_i \To +\infty$, that is,

\begin{equation}\label{Gradfi}
\int _{D(s_i)\setminus D(s_0)} \norm{\nabla f_i}^2 \To 0 , \, s_i \To + \infty.
\end{equation}

Let us denote $K^{-} = {\rm min}\set{0,K}$ and $K^+ = {\rm max}\set{0,K}$, so that
$K = K^{-} + K^{+}$. Consider the sequence of monotone functions given by
$$ g^{+}_i = K ^{+} f _i ^2 , \quad g^{-}_i = K^{-} f _i^2 $$and note that for $i=1$,
$g^{+}_1$ and $g^{-}_1 $ are integrable on $\Sigma$.


Now, apply the non positivity of $L$ to the sequence $\set{f_i}$, i.e.,
\begin{equation}\label{Kineq}
-a \int _{\Sigma } K f _i ^2 \leq \int _{\Sigma } \norm{\nabla f_i}^2 = \int
_{D(s_i)\setminus D(s_0)} \norm{\nabla f_i}^2 .
\end{equation}

We write the left hand side of this inequality as
$$ -\int _{\Sigma } K f_i ^2 = - \int _{\Sigma} g_i ^{+} - \int _{\Sigma} g_i ^{-}
.$$

By the Monotone Convergence Theorem for the sequences $\set{g_i ^{+}}$ and $\set{g_i
^{-}}$ (note that the limits could be infinite), we have
\begin{equation*}
\begin{split}
- \lim _{i \To + \infty }\int _{\Sigma } K f_i ^2 &= - \lim _{i \To + \infty
}\int _{\Sigma } g_i ^{+} -  \lim _{i \To + \infty }\int _{\Sigma } g_i ^{-} \\
 & = -\int _{\Sigma} K^{+} - \int _{\Sigma} K^{-}
\end{split}
\end{equation*}since $g_i ^± \To K^±$ uniformly on compact sets.

Thus, combining this with \eqref{Gradfi} and taking limits in \eqref{Kineq}, we have
$$ -\int _{\Sigma } K^{+} - \int _{\Sigma} K^{-} \leq 0 $$but using that $\int _{\Sigma}
K^{+}$ is finite we obtain that $\int_{\Sigma } K^{-}$ is finite, thus
\begin{equation*}
-\int _{\Sigma} K \leq 0 ,
\end{equation*}since $a$ is a positive constant.

Now, since $\int_{\Sigma} K^{-}$ is finite, the Cohn-Vosen inequality says
$$ \int _{\Sigma } K \leq 2\pi \chi (\Sigma) ,$$which means that
$$ 0 \leq \int _{\Sigma } K \leq 2 \pi \chi (\Sigma) , $$i.e., $K$ is integrable.
\end{proof}

\section{Applications to stable surfaces}

In this Section we study a non positive differential operator of the form
$$ L = \Delta + V - a K $$where $V$ is a non negative function on $\Sigma$.

If $L f \leq 0$, then the quadratic form, $I(f)$ associated to $L$ is non negative
on compactly supported functions, i.e. $I(f) \geq 0$. So, in this case, Lemma
\ref{l2} can be rewritten as

\begin{corollary}\label{c2}
Let $\Sigma $ be a Riemannian surface possibly with boundary and $K \not \equiv 0$.
Fix a point $p_0\in \Sigma$ and a positive number $ s > 0$ such that
$\overline{D(s)}\cap \partial \Sigma = \emptyset$. Suppose that the differential
operator $L = \Delta + V - a K$ is non positive on $ f\in H^{1,2}_0 (\Sigma)$, where
$V \geq 0$ and $a$ is a positive constant. Given $b\geq 1$, let $\alpha $ be defined
by \eqref{alpha}. With the notation of Lemma \ref{l2}, if $\alpha > 0$, then
\begin{equation}\label{estVpos}
\int _{D(s e^{-s})}  V  + \int _{D(s)\setminus D(s e^{-s})} \left(
\frac{\ln(s/r)}{s}\right) ^{2b} V  \leq
 2a \left(G (s) \pi + b\frac{2\pi - \mathcal{K}(se^{-s})}{s}\right) + \rho^{+}_{a,b}(\delta ,s)
\end{equation}
\end{corollary}

When $a>1/4$, we already know quadratic area growth and the integrability of the
potential (see \cite{Ca}, \cite{MPR} for $a>1/4$ or \cite{CM}, \cite{R} for
$a>1/2$).

Another interesting consequence is that we are able to bound the distance of any
point to the boundary, this is known as the \emph{Distance Lemma} (see \cite{R2},
\cite{RR} or \cite{MPR}, and \cite{Ma} for a sharp bound in space forms).

Here, we will extend this result for $0< a \leq 1/4$

\begin{theorem}\label{t5}
Let $\Sigma $ be a Riemannian surface with $k-AAG$ and possibly with boundary.
Suppose that $L = \Delta + V - a K$ is non-positive, where $V \geq c >0 $ and $0< a
\leq 1/4$. Then, there exists a positive constant $C$ such that
$$ {\rm dist}_{\Sigma}(p , \partial \Sigma) \leq C , \, \, \forall p \in \Sigma .$$

In particular, if $\Sigma$ is complete with $\partial \Sigma = \emptyset$ then it
must be topologically a sphere.
\end{theorem}
\begin{proof}
Let us suppose that the distance to the boundary were not bounded. Then there exists
a sequence of points $\set{p_i} \in \Sigma$ such that ${\rm dist}^{\Sigma}(p_i ,
\partial \Sigma) \To + \infty$. So, for each $p_i$ we can choose a real number $s_i$
such that $s_i \To + \infty$ and $\overline{D(p_i ,s_i)}\cap \partial \Sigma =
\emptyset $.

Let $\beta \in \r $ be a real number greater that one, then
$$ \int _{D(s e^{-s})}  V  + \int _{D(s)\setminus D(s e^{-s})} \left(
\frac{\ln(s/r)}{s}\right) ^{2b} V  \geq c \frac{\beta ^{2b}}{s^{2b}} a(s
e^{-\beta}). $$

Now, choose $b > 1$ such that $2(b+1) \geq k > 2b > 2 $. Thus, by \eqref{estVpos}
and the above inequality

\begin{equation}\label{eq1}
c \frac{\beta ^{2b}}{s^{2b}} a(s e^{-\beta}) \leq C  +\rho^{+}_{a,b}(\delta _0 ,s),
\end{equation}where $\delta _0$ is given in Remark \ref{rdelta0}.

Now, since $\Sigma $ has $k-$AAG and $k>2b$ then for $s$ large enough we have
\begin{equation}\label{eq2}
c \frac{\beta ^{2b}}{s^{2b}} a(s e^{-\beta}) \sim s^{k-2b} \To + \infty
\end{equation}and from \eqref{asymppos}
\begin{equation}\label{eq3}
\rho^{+}_{a ,b}(\delta_0 , s) \To \left\{ \begin{matrix}
0 & \text{for} & 2(b+1)> k \\
C^+ \rho_{{\rm min}} & \text{for} & 2(b+1)= k
\end{matrix}\right.
\end{equation}

Thus, applying \eqref{eq1} to each disk $D(p_i ,s_i )$, and bearing in mind that
from \eqref{eq2} the left hand side of \eqref{eq1} goes to infinity, and (from
\eqref{eq3}) the right hand side remains bounded, we obtain a contradiction.

We still have to consider the case $k\leq 2$. Here, we consider a formula developed
by Meeks-P\'{e}rez-Ros, this formula follows from Lemma \ref{l1} with the test
function $f(r) = \left( 1 - r/s\right)^b$ for $r \in [0,s]$, that is, for $b \geq 1$
\begin{equation}\label{MPR3}
\int _{D(s)} \left( 1- r/s\right)^{2b} V \leq 2a \pi +\frac{b(b(1-4a)+2a)}{s^2}\int
_{0}^s \left( 1- r/s\right)^{2b-2} l(r) .
\end{equation}

Thus, for $b=1$ and the $k-$AAG, $k\leq 2$, of $\Sigma$, the right hand side of
\eqref{MPR3} goes to some positive constant as $s$ goes to infinity.

But, since $V \geq c >0$,
$$ \int _{D(s)} \left( 1- r/s\right)^{2b} V \geq c \frac{a(s /2)}{4} . $$

Thus, applying \eqref{MPR3} to each disk $D(p_i ,s_i )$, and bearing in mind that
the left hand side of \eqref{MPR3} goes to infinity and the right hand side remains
bounded, we obtain a contradiction.

Note that we must be careful with the term $\mathcal{K}(s e^{-s})$ in
\eqref{estVpos}. Let us see that we do not need to worry about this term.

Let $p \in \Sigma $ be any point in the surface and consider the radial function
$u(r) = 1-r/s$ defined on $D(p ,s)$. Then, applying the non positivity of the
operator $L= \Delta + V - a K $, we have
$$ - a \int _{D(p,s)} u(r)^2  K \leq \int _{D(p,s)} \norm{\nabla u}^2  \leq \frac{a(s)}{s^2} .$$

Now, if $p$ is a point where $K(p)<0$ (note that we do not have to worry about
points where the curvature is positive), we can choose $s>0 $ small enough such that
$K(q)<0$ for all $q \in D(p, s)$, thus in this geodesic disk

$$- \frac{a}{4} \int _{D(p,s/2)}  K \leq -a \int _{D(p,s)} u(r)^2 K \leq
\frac{a(s)}{s^2},$$that is

$$ - \int _{D(p,s/2)}  K \leq \frac{4 a(s)}{a s^2} .$$

Taking into account that for $s$ small the area of $a(s) $ is almost Euclidean, we
have

$$ \frac{4 a(s)}{a s^2} \To C , \, s\To 0 $$for some positive constant $C$. So, we
finally obtain that

$$ - \int _{D(p,s/2)}  K \To C , \, s \To 0 ,$$which means that
$$ \frac{\mathcal{K}(s e^{-s})}{s} \To 0 , \, s \To + \infty . $$

Now, if $\Sigma$ is complete, then the estimate and the Hopf-Rinow Theorem imply
that $\Sigma $ must be compact. Moreover, applying the operator $L_a$ to the test
function $1$, we have
$$ a \int _{\Sigma} K \geq c \, {\rm Area}(\Sigma) $$which implies, by the Gauss-Bonnet
Theorem, that $\chi(\Sigma)>0$.
\end{proof}

As a consequence of this proof we have the following result.

\begin{corollary}
Under the hypothesis of Theorem \ref{t5}, if $-L$ has finite index, then the
distance of any point to the boundary is bounded. So, if the surface is complete, it
must be compact and its Euler characteristic is positive.
\end{corollary}

Note that following the above method we can prove (for $a > 1/4$ it is known)

\begin{theorem}\label{t6}
Let $\Sigma $ be a complete Riemannian surface satisfying $k-$AAG, $k \leq 2$.
Suppose that $L = \Delta + V - a K$ is non-positive, where $V \geq 0 $ and $0< a
\leq 1/4$. Then, $V \in L^1(\Sigma)$, i.e., $V$ is integrable.

Moreover, if $\Sigma$ has $k-$AAG with $k>2$, then for $2(b+1) \geq  k$ we have
\begin{equation*}
\int _{D(s)} V \leq C s^{2b}
\end{equation*}for some positive constant $C$.
\end{theorem}
\begin{proof}
The case when $\Sigma $ has $k-$AAG, $k\leq 2$, follows from formula \eqref{MPR3},
since then the right hand side goes to some constant, and we can bound the left hand
side as
$$ \frac{1}{4} \int _{D(s/2)} V \leq \int _{D(s)} \left( 1-r/s\right)^2 V .$$

The second case follows using that
$$ \frac{(\ln 2)^{2b}}{s^{2b}}\int _{D(s/2)} V \leq \int _{D(s e^{-s})}  V  +
\int _{D(s)\setminus D(s e^{-s})} \left( \frac{\ln(s/r)}{s}\right) ^{2b} V
$$and formula \eqref{asymppos}. So, putting this together with \eqref{estVpos}
we obtain the result.
\end{proof}

\begin{remark}
Actually, the case $k\leq 2$ in the above result has been proven in \cite{MPR} as
well.
\end{remark}

In \cite[Theorem 3]{FCS}, they proved: \emph{Let $N$ be a complete oriented
$3-$manifold of non-negative scalar curvature. Let $\Sigma$ be an oriented complete
stable minimal surface in $N$. If $\Sigma$ is noncompact, conformally equivalent to
the cylinder and the absolute total curvature of $\Sigma$ is finite, then $\Sigma$
is flat and totally geodesic.}

And they state \cite[Remark 2]{FCS}: \emph{We feel that the assumption of finite
total curvature should not be essential in proving that the cylinder is flat and
totally geodesic.}

So, using Theorem \ref{t8}, we are able to partially answer this question.

\begin{theorem}\label{t9}
Let $N$ be a complete oriented $3-$manifold of non-negative scalar curvature. Let
$\Sigma$ be an oriented complete stable minimal surface in $N$. If $\Sigma$ is
noncompact, conformally equivalent to the cylinder and $\int _{\Sigma} K^{+}$ is
finite, then $\Sigma$ is flat and totally geodesic.
\end{theorem}

\vspace{.5cm}
\begin{center}
{\bf Acknowledgement:}
\end{center}
The authors wish to thank J. P\'{e}rez and R. Schoen for their interesting comments
and help during the preparation of this work.

\end{document}